\documentclass[12pt,a4paper]{article}

\usepackage{amsthm,amsmath,natbib,amssymb,amsfonts,bm}
\usepackage{epsfig}
\usepackage{placeins}
\usepackage{epsf}
\RequirePackage[dvips]{hyperref}

\newtheorem*{theorem*}{Theorem}

\theoremstyle{definition}

\theoremstyle{remark}

\begin{document}

\baselineskip=21pt

\begin{center}
 {\bf \Large Confidence intervals for the normal mean utilizing prior information}
\end{center}

\bigskip

\begin{center}
{\bf \large David Farchione and Paul Kabaila$^{\textstyle{^*}}$}
\end{center}

\medskip

\noindent{\it $\textstyle{^*}$Department of Mathematics and
Statistics, La Trobe University, Victoria 3086, \newline Australia}

\bigskip
\noindent{\bf Abstract}
\medskip

Consider $X_{1},X_{2},\ldots,X_{n}$ that are independent and
identically $\text{N}(\mu,\sigma^{2})$ distributed. Suppose that
we have uncertain prior information that $\mu = 0$. We answer the
question: to what extent can a frequentist $1-\alpha$ confidence interval
for $\mu$ utilize this prior information?

\bigskip
\bigskip


\vbox{\vskip 6.5cm}

\noindent $^*$ Corresponding author. Address: Department of
Mathematics and Statistics, La Trobe University, Victoria 3086,
Australia; Tel.: +61-3-9479-2594; fax: +61-3-9479-2466. \newline
{\it E-mail address:} P.Kabaila@latrobe.edu.au.

\vfil \eject

\baselineskip=20pt

\noindent {\bf 1. Introduction}

\medskip

Suppose
that $X_1, \ldots, X_n$ are independent and identically $N(\mu,
\sigma^2)$ distributed. The parameter of interest is $\mu$. Also
suppose that previous experience with similar data sets and/or
scientific background and expert opinion suggest that $\mu=a$,
where $a$ is a specified number.
Without loss of generality we
assume that $a=0$.
Our aim is to answer the following question. To what extent
can a frequentist $1-\alpha$ confidence interval (i.e. a confidence interval
whose coverage probability has infimum $1-\alpha$) utilize this prior information?

For the sake of simplicity, we first deal with the case that
$\sigma^2$ is known. We find a confidence interval for $\mu$ by
first finding a confidence interval for $\theta = \sqrt{n} \mu /
\sigma$. Let $\bar X = n^{-1} \sum_{i=1}^n X_i$ and
$X = \bar X /(\sigma/\sqrt{n})$, so that
$X \sim N(\theta,1)$. Suppose that $I = \big[\ell(X),
u(X)\big]$ is a $1-\alpha$ confidence interval for $\theta$ i.e.
$\text{P}_{\theta} (\theta \in I) \ge 1-\alpha$ for all
$\theta$.
The confidence interval for $\mu$ that corresponds to this
confidence interval for $\theta$ is
$J = \big [ (\sigma/\sqrt{n}) \ell(X), \,
(\sigma/\sqrt{n}) u(X) \big ]$. Pratt (1961, 1963) considers
$X \sim N(\theta,1)$ and presents confidence
intervals for $\theta$  that (a) have a pre-specified
minimum coverage $1-\alpha$ and (b) are shorter than the usual
confidence interval when $\theta=0$.
The $1-\alpha$ confidence interval for $\mu$ that has the
smallest possible expected length when $\mu=0$ is derived by Pratt (1961) and is
\begin{equation}
\label{Pratt_interval} \left [ \min \left (0, \bar X - z_{\alpha}
\frac{\sigma}{\sqrt{n}} \right ), \max \left (0, \bar X + z_{\alpha}
\frac{\sigma}{\sqrt{n}} \right ) \right ]
\end{equation}
where $z_a$ is defined by $P(Z > z_a) = a$ for $Z \sim N(0,1)$.

This
confidence interval for $\mu$ has the following
analogue for the case that $\sigma^2$ is unknown
\begin{equation}
\label{Bofinger} \left [ \min \left (0, \bar X - t_{\alpha, n-1}
\frac{S}{\sqrt{n}} \right ), \max \left (0, \bar X + t_{\alpha, n-1}
\frac{S}{\sqrt{n}} \right ) \right ]
\end{equation}
where $t_{a,m}$ is defined by $P(T > t_{a,m}) = a$ for $T \sim
t_{a,m}$ and
$S^2 = (n-1)^{-1} \sum_{i=1}^n (X_i - \bar X)^2$.
This analogue is given by Brown et al (1995) and has been described e.g. by Bofinger (1985).

\medskip

These confidence intervals have two major problems. The first
problem is that the expected lengths of these confidence intervals
diverge to $\infty$ as $|\mu| \rightarrow \infty$. This unpleasant
feature
means that if the prior information happens to be badly incorrect
(i.e. $\mu$ happens to be very far from 0) then these confidence
intervals perform extremely poorly. The second problem is that
neither of these confidence intervals approaches the corresponding
standard confidence interval when the data strongly indicates that
the prior information about $\mu$ is incorrect. Surely, if the data
strongly indicate that this prior information is incorrect then we
should be using something very close to the standard $1 - \alpha$ confidence
interval for $\mu$. For example, when $\sigma^2$ is known and $|X| >
10$ then we should be using the standard confidence interval
$\big[\bar X - z_{\alpha/2} (\sigma/\sqrt{n}), \bar X +
z_{\alpha/2} (\sigma/\sqrt{n}) \big]$ for $\mu$.

In this paper we describe confidence intervals for $\mu$ that do not
suffer from these problems. Similarly to Hodges and Lehmann (1952)
and Bickel (1983, 1984),
our aim is to utilize the uncertain prior
information in the frequentist inference of interest,
whilst providing a safeguard in case this prior information happens to be incorrect.
Our $1 - \alpha$ confidence intervals
have the following desirable properties. They have expected lengths
that (a) are relatively small when the prior information that
$\mu=0$ is correct and (b) have a maximum value that is not too
large. They also coincide with the corresponding standard $1 - \alpha$ confidence
interval when the data happens to strongly contradict the prior
information about $\mu$. In Sections 2, 3 and 4 we deal with
the case that $\sigma^2$ is known, by applying the methodology of
Pratt (1961) with a novel weight function. In Sections 5 and 6 we use
the same novel weight function,
combined with invariance and a new computationally-based approach, to deal with 
the case that $\sigma^2$ is unknown.

Finally, consider point and interval estimators utilizing uncertain prior information
in linear regression. Bickel (1984) presents a comprehensive analysis of point estimators
in a very general setting. He also
analyzes the coverage properties of some fixed-width confidence 
intervals, assuming the covariance matrix of the error vector is 
known. Tuck (2006) develops a new variable-width confidence interval
analogous to \eqref{Pratt_interval}.
The methodology described in Sections 5 and 6 of the present paper,
leading to variable-width confidence intervals,
can be extended to the linear regression context (Kabaila and Giri (2007)).

\bigskip

\noindent {\bf 2. The known variance case}

\medskip

Assume that $\sigma^2$ is known. In the introduction we defined the
random variable $X$, which has an $N(\theta,1)$ distribution, and
the $1-\alpha$ confidence intervals $I$ and $J$ for $\theta$ and
$\mu$ respectively. Observe that $P_{\mu} (\mu \in I) =
P_{\theta} (\theta \in I)$, so that the confidence intervals $I$ and
$J$ have the same minimum coverage probabilities. Furthermore,
$E_{\mu} (\text{length of } J) =  (\sigma/\sqrt{n})
E_{\theta} (\text{length of } I)$ when $\theta =
\sqrt{n} \mu/\sigma$.
In other words, the expected length of $J$ is proportional to the
expected length of $I$. We therefore focus on the case that
$X$ has an $N(\theta,1)$ distribution and we have uncertain prior
information that $\theta=0$.

Let $C(X)$ be a $1-\alpha$ confidence region for $\theta$. Define
$A(\theta)$ by
$\theta \in C(x)$ if and only if $x \in A(\theta)$.
Here $A(\theta)$ is the acceptance region for the null hypothesis
that $\theta$ is the true parameter value. Let $L(C(X))$ denote the
sum of the lengths of the intervals making up $C(X)$.
Also let our aim be to minimise the average expected length
\begin{equation}
\label{av_exp_len_general} \int E_{\theta} \big(L(C(X)) \big) d
\nu(\theta).
\end{equation}
for a specified weight function $\nu$. We use $\phi$ to denote the
$N(0,1)$ density function. As proved by Pratt (1961), the solution to
this problem is to choose the acceptance region $A(\theta)$ to
consist of those values of $x$ such that
\begin{equation*}
\frac{\displaystyle{\int_{-\infty}^{\infty} \phi(x-\theta) d
\nu(\theta)}} {\phi(x-\theta)} < c_{\alpha}(\theta)
\end{equation*}
where $c_{\alpha}(\theta)$ is chosen such that $P_{\theta}(X \in
A(\theta)) = 1-\alpha$. For some weight functions $\nu$ the average
expected length is infinite, even for the confidence interval
corresponding to the acceptance regions obtained in this way.
However, the criterion
\begin{equation*}
\int \big (E_{\theta} \big(L(C(X)) \big) - 2 z_{\alpha/2} \big ) d
\nu(\theta)
\end{equation*}
takes the (finite) value 0 when $C(X)$ is the standard $1-\alpha$
confidence interval for $\theta$. It is straightforward to show that
the minimization of this criterion leads to the same formula for
$A(\theta)$ as the (formal) minimization of
\eqref{av_exp_len_general}.

As pointed out by Pratt (1961), the standard $1-\alpha$ confidence
interval for $\theta$,
\begin{equation}
\label{standard_CI} \big[ X - z_{\alpha/2}, X + z_{\alpha/2} \big],
\end{equation}
is the $1-\alpha$ confidence interval that minimizes the average
expected length when $\nu(x)=x$ for all $x$. Also, as pointed out by Pratt (1961),
the $1-\alpha$ confidence interval \eqref{Pratt_interval}
for $\theta$ is the $1-\alpha$ confidence interval that minimizes
the average expected length when $\nu = H$ where $H$ is the unit step
function defined by $H(x)=0$ for $x < 0$ and $H(x)=1$ for $x \ge 0$.

Now consider a weight function that is a mixture of the weight
functions $\nu(x)=x$ and $\nu = H$. It is plausible that if we
minimise the average expected length using this weight function then
we will obtain a confidence interval that (a) has relatively small
expected length when $\theta=0$ and (b) overcomes the weaknesses of
Pratt's interval \eqref{Pratt_interval}. So, we consider the
$1-\alpha$ confidence interval that minimises the average expected
length when the weight function $\nu$ is
\begin{equation}
\label{mixed_wt_fn} \nu(x) = w x + H(x) \ \text{ for all } \ x.
\end{equation}
Here, $w$ is a fixed nonnegative number. We call this the `mixed
interval'.

In this case, the acceptance region $A(\theta)$ corresponding to the
confidence region $C(X)$ minimizing the average expected length
\eqref{av_exp_len_general} consists of the values of $x$ such that
\begin{equation*}
\frac{w + \phi(x)}{\phi(x-\theta)} - c_{\alpha}(\theta) < 0,
\end{equation*}
where $c_{\alpha}(\theta)$ is chosen such that $P_{\theta}(X \in
A(\theta))=1-\alpha$. Define
\begin{equation*}
g(x,c,\theta) = \frac{w + \phi(x)}{\phi(x-\theta)} - c.
\end{equation*}
Also define $B(c,\theta) = \{ x: \, g(x,c,\theta) < 0\}$. Obviously,
$c_{\alpha}(\theta)$ is the value of $c$ such that $P_{\theta}(X \in
B(c,\theta))=1-\alpha$.
To analyse the properties of the acceptance region $A(\theta)$ we
will need the following theorem.

\medskip

\noindent {\bf Theorem 2.1.} For every fixed $w>0$, $\theta \in
\mathbb{R}$ and $c>0$, the following is true. The set $B(c,\theta)$
is either (a) the empty set or (b) an interval with finite
endpoints.

\medskip

\noindent {\bf Proof.} \ Fix $w>0$, $\theta \in \mathbb{R}$ and
$c>0$. Observe that $g(x,c,\theta) \rightarrow \infty$ as $|x|
\rightarrow \infty$. The result will be proved by showing that
$\partial g(x,c,\theta)/\partial x$ is an increasing function of $x
\in \mathbb{R}$. Now
$g(x,c,\theta) = \exp(\textstyle{\frac{1}{2}} \theta^2) g^*(x,\theta)
- c$, where
$g^*(x,\theta) = w^* \exp(\textstyle{\frac{1}{2}} x^2 - \theta x) +
\exp(-\theta x)$
and $w^* = \sqrt{2 \pi} w$. Note that
\begin{equation*}
\frac{\partial g^*(x,\theta)}{\partial x} =
\exp(-\textstyle{\frac{1}{2}} \theta^2) \, w^* (x-\theta)
\exp(\textstyle{\frac{1}{2}} (x-\theta)^2) - \theta \exp(-\theta x).
\end{equation*}
This is an increasing function of $x$. Hence $\partial
g(x,c,\theta)/ \partial x$ is an increasing function of $x$.

\hfill $\square$

\noindent This theorem leads to the very important property of
$A(\theta)$ described in the following corollary, whose proof is
omitted for the sake of brevity.

\medskip

\noindent {\bf Corollary 2.1.} For every fixed $w>0$ and $\theta \in
\mathbb{R}$, the following is true. The $1-\alpha$ acceptance region
$A(\theta)$ is an interval with finite endpoints.

\medskip

\noindent The computation of the acceptance region $A(\theta)$ for
given $w>0$ and $\theta \in \mathbb{R}$ is facilitated by the
following result.

\medskip

\noindent {\bf Theorem 2.2.} For every $w>0$ and $\theta \in
\mathbb{R}$,
\begin{equation*}
w \sqrt{2 \pi} \exp(\textstyle{\frac{1}{2}} z^2_{\alpha/2}) \le
c_{\alpha}(\theta) \le (w \sqrt{2 \pi} + 1)
\exp(\textstyle{\frac{1}{2}} z^2_{\alpha/2}).
\end{equation*}

\medskip

\noindent {\bf Proof.} The result follows from the fact that for
every $w>0$ and $\theta \in \mathbb{R}$ the following is true. For
every $x \in \mathbb{R}$,
\begin{equation*}
\frac{w}{\phi(x-\theta)} \le \frac{w + \phi(x)}{\phi(x-\theta)} \le
\frac{w + (1/\sqrt{2 \pi})}{\phi(x-\theta)}.
\end{equation*}
\hfill $\square$

\noindent The following theorem describes an important property of the 
confidence set $C(x)$. The proof of this theorem is omitted, for the
sake of brevity.

\medskip

\noindent {\bf Theorem 2.3.} For every $w > 0$ the following is
true. The $1-\alpha$ confidence set $C(x)$ is an interval for all
sufficiently large $|x|$, with endpoints approaching those of the
standard $1-\alpha$ confidence interval $\big[x - z_{\alpha/2},x +
z_{\alpha/2} \big]$ as $|x| \rightarrow \infty$.


\baselineskip=20pt
\bigskip

\noindent {\bf 3. Numerical comparison of the intervals for the
known variance case}

\medskip

We continue with our consideration of the case that $\sigma^2$ is
known. As described in the introduction, we reduce this case to the
problem of finding a $1-\alpha$ confidence interval for $\theta$
based on $X \sim N(\theta,1)$. We denote the standard $1-\alpha$
confidence interval \eqref{standard_CI} by $C_S(X)$. Also, we denote
Pratt's interval \eqref{Pratt_interval} by $C_P(X)$.

Consider the case that the weight function $\nu$ is given by
\eqref{mixed_wt_fn}. This weight function is a mixture of the weight
functions $\nu(x)=x$ and $\nu = H$ that lead to $C_S$ and $C_P$
respectively. For $1-\alpha = 0.95$ and each $w \in \{ 0.01, 0.1, 1
\}$, Corollary 2.1 and Theorem 2.2 were used to compute the
acceptance regions $A(\theta)$, corresponding to the 0.95 confidence
sets minimizing the average expected length
\eqref{av_exp_len_general}, for a fine grid of values of $\theta$.
For each of these values of $w$, the confidence region corresponding
to $A(\theta)$ was found to always be an interval.
We denote the 0.95 confidence interval minimizing the
average expected length when $\nu$ is given by \eqref{mixed_wt_fn}
(which we have called the mixed interval) by $C_M^w(X)$. All the
computations for the present paper were performed with programs
written in MATLAB, using the Optimization and Statistics toolboxes.

We use $C_S$ as the standard against with other $1-\alpha$
confidence intervals for $\theta$ are judged. The efficiency of
$C_S$ relative to $C$, a $1-\alpha$ confidence interval for
$\theta$, for a given value of $\theta$ is defined to be
\begin{equation*}
e(\theta, C_S, C) = \left (\frac{E_{\theta} \big(L(C(X))
\big)}{E_{\theta} \big(L(C_S(X)) \big)} \right )^2.
\end{equation*}
Let $\bm{X} = (X_{1},X_{2}, \ldots ,X_{n})$. Note that a $1-\alpha$
confidence interval $C(X)$ for $\theta$ corresponds to a $1-\alpha$
confidence interval $D(\bm{X})$ for $\mu$ that is obtained by
multiplying the endpoints of $C(X)$ by $\sigma/\sqrt{n}$. Let
$D_S(\bm{X})$ denote the standard $1-\alpha$ confidence interval for
$\mu$. We define the efficiency of $D_S$ relative to $D$ as
$\big (E (L(D(\bm{X}))) / E (L(D_S(\bm{X})) )
\big )^2$
and note that this is the same function of $\theta$ as $e(\theta,
C_S, C)$.

Figure 1 shows plots of the efficiency of $C_S$ relative to $C_M^w$
for $w=1$, $w=0.1$, $w=0.01$ and $w=0$, when $1-\alpha=0.95$. Note
that Pratt's interval $C_P$ is equal to the mixed interval $C_M^w$
for $w=0$. Also, the standard interval $C_S$ may be viewed as the
mixed interval $C_M^w$ in the limiting case $w \rightarrow \infty$.
Even for $w=1$, which is not a particularly large value of $w$,
$C_M^w$ is fairly close to $C_S$.

The minimum over all $1-\alpha$ confidence intervals $C$ of $e(0,
C_S, C)$ is 0.7223 and this minimum is achieved by Pratt's interval
$C_P$. However, as noted in the introduction, this interval suffers
the severe problems that (a) $e(\theta, C_S, C_P) \rightarrow
\infty$ as $|\theta| \rightarrow \infty$ and (b) it does not revert
to the standard interval when $|X| \rightarrow \infty$. The mixed
interval $C_M^w$ with $w=0.1$ is far preferable. For this interval,
$e(0, C_S, C_M^w) = 0.8016$, which is not that much larger than
0.7223. Also, for this interval, $e(\theta, C_S, C_M^w)$ never
exceeds 1.2095. Finally, in accordance with Theorem 2.3, this
interval approaches the standard interval $C_S$ as $|X| \rightarrow
\infty$. Of course, the value of $w$ can be chosen so as to reflect
the strength of the prior information that $\theta=0$.

\begin{figure}[h]
\centering \hskip-3mm
\includegraphics[scale=1.0]{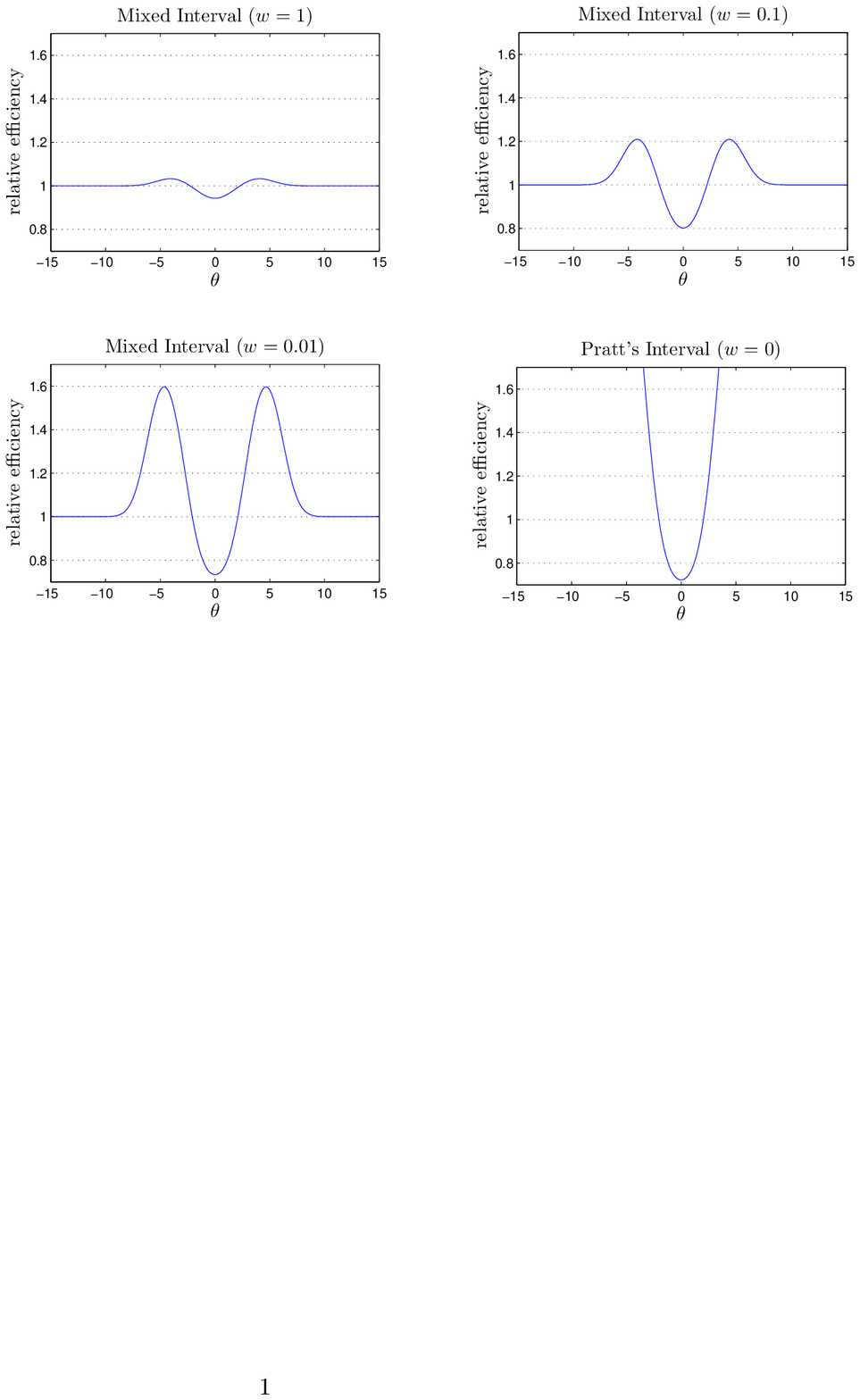} \caption{Plots
of the efficiency $e(\theta, C_S, C_M^w)$ of the standard interval
$C_S$ relative to the mixed interval $C_M^w$ for $w=1$, $w=0.1$,
$w=0.01$ and $w=0$ when $1-\alpha=0.95$.}\label{Fig2}
\end{figure}

\FloatBarrier

\baselineskip=20pt

\noindent {\bf 4. Invariance properties of the confidence interval
in the known
variance case}

\medskip

In this section we first describe the invariance properties that we
expect the confidence interval $J$ (defined in the introduction)
to possess. Traditional invariance
arguments (see e.g. Casella and Berger (2002, section 6.4) do not include
considerations of the available prior information. The novelty in
the present section is that the invariance arguments need to take
proper account of the uncertain prior information that $\mu = 0$.

We first describe an invariance property that $J$ already possesses.
The model that $X_1,\ldots, X_n$ are independent and identically
$N(\mu, \sigma^2)$ distributed may be re-expressed as follows.
Define $Y_i = a X_i$ for $i=1,\ldots, n$ where $a >0$. Thus $Y_1,
\ldots, Y_n$ are independent and identically $N(\tilde \mu, \tilde
\sigma^2)$ distributed where $\tilde \mu = a \mu$ and $\tilde \sigma
= a \sigma$. Define $Y = \bar Y / (\tilde \sigma/\sqrt{n})$. The
prior information $\mu=0$ may be re-expressed as
$\tilde \mu=0$. The re-expressed model and
prior information have the same form as the original model and prior
information. Thus the confidence interval
$\tilde J = \big [ (\tilde \sigma/\sqrt{n}) \ell(Y), \, (\tilde \sigma/\sqrt{n}) u(Y) \big ]$
for $\tilde \mu$ must lead to a confidence interval for $\mu$ that
is identical to $J$. It is easily seen that this requirement is
satisfied.

Next, we describe an invariance property that $J$ should possess and
conclude from this that the equality $\ell(x) = - u(-x)$ should hold
for all $x \in \mathbb{R}$. The model that $X_1,\ldots, X_n$ are
independent and identically $N(\mu, \sigma^2)$ distributed may be
re-expressed as follows. Define $Y_i = - X_i$ for $i=1,\ldots, n$.
Thus $Y_1, \ldots, Y_n$ are independent and identically $N(\tilde
\mu, \sigma^2)$ distributed where $\tilde \mu = - \mu$. Define $Y =
\bar Y / (\sigma/\sqrt{n})$. The prior information
$\mu=0$ may be re-expressed as $\tilde \mu=0$. The re-expressed model and prior information have the same
form as the original model and prior information. Thus the
confidence interval
$\tilde J = \big [ (\sigma/\sqrt{n}) \ell(Y), \,
(\sigma/\sqrt{n}) u(Y) \big ]$
for $\tilde \mu$ must lead to a confidence interval for $\mu$ that
is identical to $J$. It is easily seen that this requirement is
satisfied if and only if $\ell(x) = - u(-x)$ for all $x \in
\mathbb{R}$. Note that both Pratt's interval \eqref{Pratt_interval}
and the mixed interval defined in Section 2 satisfy this
requirement.

\bigskip

\noindent {\bf 5. Invariance arguments in the unknown variance case}

\medskip

We now consider the case that $\sigma^2$ is unknown. This is the
case that is usually encountered in practice. The standard
$1-\alpha$ confidence interval for $\mu$ is
\begin{equation}
\label{standard_CI_unknown_var} \left [ \bar X - t_{\alpha/2, n-1}
\frac{S}{\sqrt{n}}, \bar X + t_{\alpha/2, n-1} \frac{S}{\sqrt{n}}
\right ].
\end{equation}
A natural analogue of the confidence interval $J$ for $\mu$ is the
confidence interval
\begin{equation*}
K = \left [\frac{S}{\sqrt{n}} a\left
(\frac{\bar X}{S/\sqrt{n}}\right), \frac{S}{\sqrt{n}} b\left
(\frac{\bar X}{S/\sqrt{n}}\right) \right ]
\end{equation*}
for $\mu$. Note that both \eqref{standard_CI_unknown_var} and
\eqref{Bofinger} have this form. Suppose that our uncertain prior
information is that $\mu=0$. Using the same model transformations as in Section 4,
invariance arguments show that the equality $a(x) = -b(-x)$ must hold for
all $x \in \mathbb{R}$. In other words,
\begin{equation}
\label{final_K} K = \left [-\frac{S}{\sqrt{n}} b\left (\frac{-\bar
X}{S/\sqrt{n}}\right), \frac{S}{\sqrt{n}} b\left (\frac{\bar
X}{S/\sqrt{n}}\right) \right ].
\end{equation}
The constraint that the upper endpoint of this confidence interval
is never less than the lower endpoint implies that $b(x) \ge -b(-x)$
for all $x \in \mathbb{R}$. It also seems reasonable to require that
$b$ is a strictly increasing function.

\bigskip

\noindent {\bf 6. Computation of the interval for the unknown
variance case}

\medskip

In this section we provide computationally convenient expressions
that are used to calculate the mixed interval for the unknown
variance case. 
We illustrate the performance of this interval and compare its performance
with the
corresponding mixed interval when $\sigma^2$ is known for
the case that $n=24$ and $1-\alpha=0.95$.

Suppose that $\sigma^{2}$ is unknown. As in Section 3, let $\bm{X} =
(X_{1},X_{2}, \ldots ,X_{n})$. Also let $G(\bm{X})$ be a  confidence
interval for $\mu$ that is of the form \eqref{final_K}. Our aim is
to minimize the average expected length of $G(\bm{X})$ for a given
weight function $\nu$, such that the coverage probability of
$G(\bm{X})$ is at least $1 - \alpha$ for all $\mu \in \mathbb{R}$.
Let $\theta = \sqrt{n}\mu/\sigma$. As we show shortly, the coverage
probability $P(\mu \in G(\bm{X}))$ is a function of $\theta$. The
expected length of $G(\bm{X})$ is a function of $(\mu, \sigma^2)$.
However, we will introduce a scaled expected length of $G(\bm{X})$
which is a function of $\theta$. By using this scaled expected
length, instead of the expected length, we will be able to achieve
our aim by considering only quantities that are functions of
$\theta$ (cf. Kabaila (1998, 2005)).

The coverage probability of $G(\bm{X})$ is a function of $\theta$
and is given by
\begin{equation}
\label{Coverage Probability} P(\mu \in G(\bm{X})) = P \left(-R \, b
\left ( \frac{-X}{R} \right ) \le \theta \le R \, b \left (
\frac{X}{R} \right ) \right )
\end{equation}
where $X = \sqrt{n}\bar{X}/\sigma \sim \text{N}(\theta,1)$ and $R =
S/\sigma$.  Note that $X$ and $R$ are independent random variables.
We assume that $b$ is a strictly increasing function. This implies
that $b^{-1}$ exists. A computationally convenient expression for the
right hand side of \eqref{Coverage Probability} is
\begin{equation}
\label{Computational Coverage Probability} \int_{0}^{\infty}
\left(\Phi\left(-r b^{-1}\left(\frac{-\theta}{r}\right) -
\theta\right) - \Phi\left(r b^{-1}\left(\frac{\theta}{r}\right) -
\theta\right)\right) f_{R}(r) dr
\end{equation}
where $\Phi$ is the N(0,1) cumulative distribution function and
$f_{R}$ denotes the density of $R$.

We introduce the scaled expected length of $G(\bm{X})$ which is a
function of $\theta$ and is defined to be
\begin{equation*}
\frac{\sqrt{n}}{\sigma} \text{E}\left(L(G(\bm{X}))\right) =
\text{E}\left(R\left(b\left(\frac{X}{R}\right) +
b\left(\frac{-X}{R}\right)\right)\right).
\end{equation*}
A computationally convenient expression for this scaled expected
length is
\begin{equation}
\label{sc_exp_len_conv}
\int_{0}^{\infty}\int_{-\infty}^{\infty}
\left(b\left(\frac{x}{r}\right) + b\left(\frac{-x}{r}\right)\right)
\phi(x - \theta) \ dx \ rf_{R}(r) dr.
\end{equation}

We use the weight function \eqref{mixed_wt_fn}. As with the
$\sigma^{2}$ known case, for $w>0$, the average scaled expected
length criterion is infinite even for the standard confidence
interval \eqref{standard_CI_unknown_var}. Therefore, similarly to
Section 2, we replace this criterion by the following criterion
\begin{equation}
\label{Weighted Difference} \int \left(\frac{\sqrt{n}}{\sigma}
\text{E}\left(L(G(\bm{X}))\right) - 2 t_{\alpha/2,n - 1} \text{E}(R)
\right) d \nu(\theta)
\end{equation}
which takes the (finite) value 0 when $G(\bm{X})$ is the standard
confidence interval \eqref{standard_CI_unknown_var}. Substituting
the expression \eqref{sc_exp_len_conv} for the scaled expected
length into \eqref{Weighted Difference} we obtain
\begin{equation}
\label{Computational Weighted Difference1}
\int_{0}^{\infty}\int_{-\infty}^{\infty}
\left(b\left(\frac{x}{r}\right) + b\left(\frac{-x}{r}\right) - 2
t_{\alpha/2,n - 1}\right) (w + \phi(x)) \ dx \ rf_{R}(r) dr.
\end{equation}
Remember that we require the confidence interval to coincide with
the standard $1 - \alpha$ confidence interval when the data happens
to strongly contradict the prior information about $\mu$. The
statistic $\sqrt{n} \bar X/S$ provides an indication of how far
$\sqrt{n} \mu/\sigma$ is from 0. We therefore satisfy this
requirement by setting $b(y) = y + t_{\alpha/2,n - 1}$ for all $|y|
\ge q$ where $q$ is a specified positive number (which is chosen to
be sufficiently large). Thus $b(x/r) + b(-x/r) - 2t_{\alpha/2,n - 1}
= 0$ for all $|x|/r \ge q$. Changing the variable of integration
from $x$ to $y = x/r$, \eqref{Computational Weighted Difference1}
can now be expressed in the computationally convenient form
\begin{equation}
\label{Computational Weighted Difference2}
\int_{0}^{\infty}\int_{-q}^{q} \big(b(y) + b(-y) - 2 t_{\alpha/2,n -
1}\big) (w + \phi(ry)) \ dy \ r^{2}f_{R}(r) dr.
\end{equation}
For computational ease, we restrict the function $b(y)$ to be a
cubic spline in the interval $[-q,q]$. This spline is required to
take the value $-q + t_{\alpha/2,n - 1}$ at $y=-q$ and $q +
t_{\alpha/2,n - 1}$ at $y=q$ and has knots that are equally spaced
between $-q$ and $q$. In addition, the derivative of this spline is
constrained to be 1 at both $y=-q$ and $y=q$.

We minimize \eqref{Computational Weighted Difference2} with respect
to the function $b$, subject to the constraints on $b$ described at
the end of Section 5 and the constraint that \eqref{Computational
Coverage Probability} is at least $1-\alpha$ for all $\theta \in
\mathbb{R}$. We denote the minimizing function $b$ by $b_{M}^{w}$.
We call the confidence interval for $\mu$ corresponding to
$b_{M}^{w}$ the mixed interval and denote it by $G_{M}^{w}$. We
denote the standard confidence interval
\eqref{standard_CI_unknown_var} by  $G_{S}$.  Similarly to
Section 3, we use $G_S$ as the standard against with other
$1-\alpha$ confidence intervals for $\mu$ are judged. The efficiency
of $G_S$ relative to $G$, a $1-\alpha$ confidence interval for
$\mu$, is defined to be
$\left (E \big(L(G(\bm{X})) \big)/ E \big(L(G_S(\bm{X})) \big) \right )^2$
which is a function of $\theta$.

\begin{figure}[b]
\centering \hskip-5mm
\includegraphics[scale=1.1]{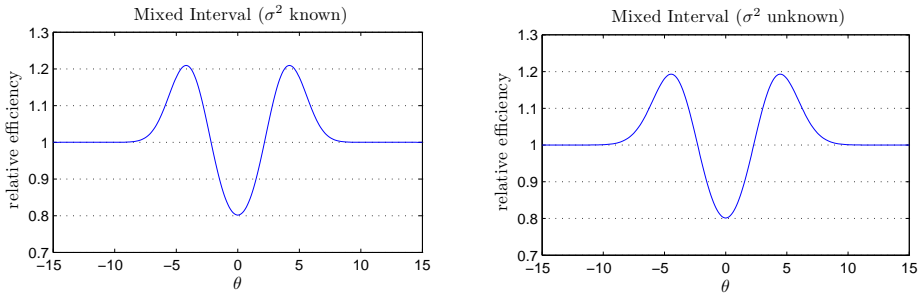} \vskip-5mm\caption{{These plots
concern the case that $n = 24$, $1-\alpha=0.95$ and $w = 0.1$. The
left panel is a plot of the efficiency of $D_S$ relative to
$D_M^w$ as a function of $\theta$. The right panel is a plot of
the efficiency of $G_S$ relative to $G_M^w$ as a function of
$\theta$. }}\label{Fig2}
\end{figure}

To illustrate the performance of the mixed interval and to compare its performance
with the
corresponding mixed interval when $\sigma^2$ is known, we consider
the case that $n=24$ and $1-\alpha=0.95$.
For the computation of $G_M^w$, we chose $q=8$ with the knots of the
cubic spline at $-8, -7, \ldots, 7, 8$. We also chose $w = 0.1$. The
efficiency of $G_S$ relative to $G_M^w$ is shown in the right
panel of Figure 2. When the prior information is correct i.e.
$\mu=0$ the efficiency of $G_S$ relative to $G_M^w$ is 0.8013. Also,
the efficiency of $G_S$ relative to $G_M^w$ never exceeds 1.1930.
Furthermore, $G_M^w$ reverts to the standard $1-\alpha$ confidence
interval when the prior information happens to be badly incorrect.
This is reflected in the fact that the efficiency of $G_S$ relative
to $G_M^w$ approaches 1 as $|\theta| \rightarrow \infty$. It is
notable that the coverage probability of the confidence interval
$G_M^w$ was found to be 0.95 to an extremely good approximation
throughout the parameter space.
Now $n=24$ is quite large and so $S$ will be
probabilistically close to $\sigma$. We therefore expect that 
the efficiency of $G_S$ relative to $G_M^w$ to be similar to the
efficiency of $D_S$ relative to $D_M^w$ when $w=0.1$. This expectation is
confirmed by the left panel of Figure 2.

\FloatBarrier
\baselineskip=21pt
\noindent {\bf References}

\smallskip

\noindent Bickel, P.J., (1983). Minimax estimation of the mean of a normal distribution subject
to doing well at a point, in: M.H. Rizvi, J.S. Rustagi and D. Siegmund, eds. Recent Advances in Statistics, Academic
Press, New York, pp. 511--528.

\noindent Bickel, P.J., (1984). Parametric robustness: small biases can be worthwhile.
Annals of Statistics 12, 864--879.

\noindent Bofinger, E., 1985. Expanded confidence intervals.
Communications in Statistics: Theory and Methods 14, 1849--1864.

\noindent Brown, L.D., Casella, G., Hwang, J.T.G., (1995). Optimal
confidence sets, bioequivalence and the Limacon of Pascal.
Journal of the American Statistical Association 90,880--889.

\noindent Casella, G., Berger, R. L., (2002). \textit{Statistical
Inference, 2nd ed.} Duxbury, Pacific Grove, California.

\noindent Hodges, J.L., Lehmann, E.L., (1952). The use of previous
experience in reaching statistical decisions. Annals of Mathematical
Statistics 23,396--407.

\noindent Kabaila, P., (1998). Valid confidence intervals in
regression after variable selection. Econometric Theory
14, 463--482.

\noindent Kabaila, P., (2005). Assessment of a preliminary F-test
solution to the Behrens-Fisher problem. Communications in Statistics:
Theory and Methods 34, 291--302.

\noindent Kabaila, P., Giri, K., (2007). Large sample confidence 
intervals in regression utilizing prior information. La Trobe University,
Department of Mathematics and Statistics, Technical Report No. 2007--1, Jan
2007.

\noindent Pratt, J.W., (1961). Length of confidence intervals.
Journal of the American Statistical Association 56, 549--657.

\noindent Pratt, J.W., (1963). Shorter confidence intervals for the
mean of a normal distribution with known variance. Annals of
Mathematical Statistics 34, 574--586.

\noindent Tuck, J., (2006). Confidence intervals incorporating prior
information. PhD thesis, August 2006, Department of Mathematics and
Statistics, La Trobe University.

\end{document}